\documentclass[reqno]{amsart}
\usepackage{amssymb}
\usepackage{amsfonts}

\global\baselineskip24pt plus2pt minus2pt
\global\baselineskip24pt plus2pt minus2pt

\theoremstyle{plain}

\newtheorem{definition}{Definition}

\numberwithin{equation}{section}
\input{tcilatex}

\begin{document}
\title[Hyperbolic and Elliptic Equations]{On Finite Product of Convolutions and Classifications of
Hyperbolic and Elliptic
Equations}
\author{A. K{\i}l{\i}\c cman and H. Eltayeb}
\address{Department of Mathematics and Institute for Mathematical Research,
University Putra Malaysia, 43400 UPM, Serdang, Selangor, Malaysia}
\email{akilicman@umt.edu.my,\ \ hassangadeen@hotmail.com}
\date{\today }

\begin{abstract}
In this paper we consider the linear second order partial
differential equation with non-constant coefficients; then by
using the double convolution product we produce new equations with
polynomials coefficients and we classify the new equations. It is
shown that the classifications of hyperbolic and elliptic new
equations are similar to the original equations that is the
classification is invariant after finite double convolutions
product.
\end{abstract}

\maketitle

\section{\protect\bigskip Introduction}

\noindent There is no general method to solve a PDE and the
behavior of the solutions depend essentially on the classification
of PDE therefore the problem of classification of partial
differential equations is very natural and well known. Some of
second-order linear partial differential equations can be
classified as Parabolic, Hyperbolic or Elliptic. The
classification provides a guide to appropriate initial and
boundary conditions, as well as to smoothness of the solutions. If
a PDE has coefficients which are not constant, it is possible that
it will not belong to any of these categories but rather be of
mixed type. Recently, A. K\i l\i \c{c}man and H. Eltayeb in
\cite{Adem}, introduced a new method of classification for partial
differential equations with polynomial coefficients. In this study
we extend the current classification to the finite product of
double convolution with variable coefficients. First of all we
give the following definition since it has a relation with the
present study.

\begin{definition} {\rm If $F_{1}(x,y)$ and $F_{2}(x,y)$ are integrable functions then the
following integral
\begin{equation*}
F_{1}(x,y)\ast ^{x}F_{2}(x,y)=\int_{0}^{x}F_{1}(x-\theta ,y)F_{2}(\theta
,y)d\theta
\end{equation*}%
is called a single convolution with respect to $x$ provided that
the integral exists.  \\

\noindent Similarly, if $F_{1}(x,y)$ and $F_{2}(x,y)$ are
integrable functions then the following integral
\begin{equation*}
F_{1}(\theta _{1},\theta _{2})\ast \ast F_{2}(\theta _{1},\theta
_{2})=\int_{0}^{y}\int_{0}^{x}F_{1}(x-\theta _{1},y-\theta _{2})F_{2}(\theta
_{1},\theta _{2})d\theta _{1}d\theta _{2}
\end{equation*}%
is known as double convolution provided that the integral exists,
see the details in \cite{IN}.}
\end{definition}

\noindent Let us consider the linear second order partial
differential equation \ with non-constant coefficients in the form
of
\begin{equation}
a(x,y)u_{xx}+b(x,y)u_{xy}+c(x,y)u_{yy}+d(x,y)u_{x}+e(x,y)u_{y}+f(x,y)u=0
\label{eq:tr4.1}
\end{equation}%
\noindent and almost linear equation in two variable
\begin{equation}
au_{xx}+bu_{xy}+cu_{yy}+F(x,y,u,u_{x},u_{y})=0  \label{eq:tr4.2}
\end{equation}%
\noindent where $a,b,c,$ are polynomials and defined by
\begin{equation*}
a(x,y)=\sum_{\beta =1}^{n}\sum_{\alpha =1}^{m}x^{\alpha }y^{\beta },\text{ \
\ \ }b(x,y)=\sum_{j=1}^{n}\sum_{i=1}^{m}x^{\zeta }y^{\eta }\text{, \ }%
c(x,y)=\sum_{l=1}^{n}\sum_{k=1}^{m}x^{k}y^{l}
\end{equation*}%
\noindent

and $(a,b,c)\neq (0,0,0)$ and the expression $au_{xx}+2bu_{xy}+cu_{yy}$ is
called the principal part of Eq (\ref{eq:tr4.2}), since \noindent the
principal part mainly determines the properties of solution. If we multiply
the Eq (\ref{eq:tr4.2}) by polynomial with single convolution as $k(x)\ast
^{x}$ where $\displaystyle k(x)=\sum_{i=1}^{m}x^{i}$ then Eq (\ref{eq:tr4.2}%
) becomes
\begin{equation}
k(x)\ast ^{x}\left[
a(x,y)u_{xx}+b(x,y)u_{xy}+c(x,y)u_{yy}+F(x,y,u,u_{x},u_{y})\right] =0
\label{eq:tr444}
\end{equation}%
\noindent where the symbol \ $\ast ^{x}$ \ indicates single convolution with
respect to $x$, and we shall classify more general form Eq (\ref{eq:tr444})
instead of Eq (\ref{eq:tr4.2}) by considering the function
\begin{equation}
D(x,y)=\left( k(x,y)\ast ^{x}b(x,y)\right) ^{2}-\left( k(x,y)\ast
^{x}a(x,y)\right) \left( k(x,y)\ast ^{x}c(x,y)\right)  \label{eq:tr445}
\end{equation}%
From Eq (\ref{eq:tr445}), one can see that if $D$ is positive then Eq (\ref%
{eq:tr445}) it called Hyperbolic, if $D$ is negative then Eq (\ref{eq:tr445}%
) it called Elliptic, otherwise parabolic. \newline
\noindent First of all, we compute the coefficients of Eq (\ref{eq:tr444})
as follow
\begin{equation*}
A(x,y)=k(x)\ast ^{x}a(x,y)=\sum_{i=1}^{m}x^{i}\ast ^{x}\sum_{\beta
=1}^{n}\sum_{\alpha =1}^{m}x^{\alpha }y^{\beta }
\end{equation*}%
\noindent by using single convolution definition and integral by part we
obtain%
\begin{equation}
A(x,y)=\sum_{\beta =1}^{n}y^{\beta }\int_{0}^{x}\sum_{i=1}^{m}\sum_{\alpha
=1}^{m}(x-\mu )^{i}\mu ^{\alpha }d\mu  \label{eq:tr446}
\end{equation}%
\noindent by computing the integral in the right hand side of Eq (\ref%
{eq:tr446}) we have
\begin{equation*}
\sum_{\beta =1}^{n}y^{\beta }\int_{0}^{x}\sum_{i=1}^{m}\sum_{\alpha
=1}^{m}(x-\mu )^{i}\mu ^{\alpha }d\mu =\sum_{\beta
=1}^{n}\sum_{i=1}^{m}\sum_{\alpha =1}^{m}\frac{i!x^{\alpha +i+1}y^{\beta }}{%
(\alpha +1)((\alpha +2){...}(\alpha +i+1)}
\end{equation*}%
\noindent thus, we obtain the first coefficient of Eq (\ref{eq:tr444}) in
the form of
\begin{equation*}
A_{1}(x,y=\sum_{\beta =1}^{n}\sum_{i=1}^{m}\sum_{\alpha =1}^{m}\frac{%
i!x^{\alpha +i+1}y^{\beta }}{\left( (\alpha +1)((\alpha +2){...}(\alpha
+i+1)\right) }
\end{equation*}%
\noindent\ similarly, for\ the coefficients of the second part in Eq (\ref%
{eq:tr444}) we have
\begin{equation*}
B_{1}(x,y)=\sum_{j=1}^{n}\sum_{\zeta =1}^{m}\sum_{i=1}^{m}\frac{i!x^{\zeta
+i+1}y^{\eta }}{\left( (\zeta +1)((\zeta +2){...}(\zeta +i+1)\right) }
\end{equation*}%
\noindent also the last coefficient of Eq (\ref{eq:tr444}) given by
\begin{equation*}
C_{1}(x,y)=\sum_{l=1}^{n}\sum_{k=1}^{m}\sum_{i=1}^{m}\frac{i!x^{k+i+1}y^{l}}{%
\left( (k+1)((k+2){...}(k+i+1)\right) }
\end{equation*}%
then one can easily set up
\begin{equation}
D_{1}(x,y)=B_{1}^{2}-A_{1}(x,y)C_{1}(x,y).  \label{eq:tr447}
\end{equation}

for example in particular we can have
\begin{equation}
x^{3}\ast ^{x}x^{2}y^{3}u_{xx}+x^{3}\ast ^{x}x^{3}y^{4}u_{xy}+x^{3}\ast
^{x}x^{4}y^{5}u_{yy}=f(x,y)\ast ^{x}g(x,y).  \label{eq:tr666}
\end{equation}%
The first coefficients of Eq (\ref{eq:tr666}) given by
\begin{equation}
A_{1}(x,y)=x^{3}\ast ^{x}x^{2}y^{3}=y^{5}\int_{0}^{x}(x-\theta )^{3}\theta
^{2}d\theta =\frac{1}{60}y^{3}x^{6}.  \label{eq:tr661}
\end{equation}

Similarly, the second coefficient given by%
\begin{equation}
B_{1}(x,y)=x^{3}\ast ^{x}x^{3}y^{4}=\frac{1}{140}y^{4}x^{7}.
\label{eq:tr662}
\end{equation}

By the same way we get the last coefficients of Eq (\ref{eq:tr666})%
\begin{equation}
C_{1}(x,y)=x^{3}\ast ^{x}x^{4}y^{5}=\frac{1}{280}y^{5}x^{8}.
\label{eq:tr663}
\end{equation}%
By using Eqs (\ref{eq:tr447}), (\ref{eq:tr661}), (\ref{eq:tr662}) and (\ref%
{eq:tr663}) we obtain
\begin{equation}
D_{1}(x,y)=-\frac{1}{117600}y^{8}x^{14}  \label{eq:tr664}
\end{equation}

We can easy see from Eq(\ref{eq:tr664}) that Eq (\ref{eq:tr666}) is an
elliptic equation for all $(x_{0,}y).$\newline

In the same way, if we multiply the Eq (\ref{eq:tr4.2}) by polynomial with
single convolution as $h(y)\ast ^{y}$ where $\displaystyle %
h(y)=\sum_{j=1}^{n}y^{j}$ then Eq (\ref{eq:tr4.2}) becomes%
\begin{equation}
h(y)\ast ^{y}\left[
a(x,y)u_{xx}+b(x,y)u_{xy}+c(x,y)u_{yy}+F(x,y,u,u_{x},u_{y})\right] =0
\label{eq:tr770}
\end{equation}%
where the symbol \ $\ast ^{y}$\ indicates single convolution with respect to
$y$, and we shall classify Eq (\ref{eq:tr770}) as. First of all, let us
compute the coefficients of Eq (\ref{eq:tr770}) by using definition of
single convolution with respect to $y$ and integral by part we obtain the
first coefficient of Eq (\ref{eq:tr770}) as follow
\begin{equation*}
A_{2}(x,y=h(y)\ast ^{y}a(x,y)=\sum_{\beta =1}^{n}\sum_{j=1}^{n}\sum_{\alpha
=1}^{m}\frac{j!x^{\alpha }y^{\beta +j+1}}{\left( (\beta +1)((\beta +2){...}%
(\beta +j+1)\right) },
\end{equation*}

and the second coefficients of Eq(\ref{eq:tr770}) given by
\begin{equation*}
B_{2}(x,y)=h(y)\ast ^{y}b(x,y)=\sum_{j=1}^{n}\sum_{\eta =1}^{n}\sum_{i=1}^{m}%
\frac{j!x^{\zeta }y^{\eta +j+1}}{\left( (\eta +1)((\eta +2){...}(\eta
+i+1)\right) }.
\end{equation*}

Similarly, the last coefficients of Eq(\ref{eq:tr770}) given by
\begin{equation*}
C_{2}(x,y)=h(y)\ast ^{y}c(x,y)=\sum_{j=1}^{n}\sum_{l=1}^{n}\sum_{i=1}^{m}%
\frac{j!x^{k}y^{l+j+1}}{\left( (l+1)((l+2){...}(l+i+1)\right) }
\end{equation*}

In particular, let us classify the following example
\begin{equation}
y^{7}\ast ^{y}x^{2}y^{3}u_{xx}+y^{7}\ast ^{y}x^{3}y^{4}u_{xy}+y^{7}\ast
^{y}x^{4}y^{5}u_{yy}=f(x,y)\ast ^{y}g(x,y).  \label{eq:tr880}
\end{equation}

the symbol $\ast ^{y}$means single convolution with respect to $y$. We
follow the same technique that used above, then the first coefficient of Eq (%
\ref{eq:tr880}) given by%
\begin{equation}
A_{2}(x,y=y^{7}\ast ^{y}x^{2}y^{3}=\frac{1}{1320}x^{2}y^{11},
\label{eq:tr881}
\end{equation}

the second coefficient of (\ref{eq:tr880}) given by%
\begin{equation}
B_{2}(x,y)=y^{7}\ast ^{y}x^{3}y^{4}=\frac{1}{3860}x^{3}y^{12}
\label{eq:tr882}
\end{equation}

and the last coefficient given by
\begin{equation}
C_{2}(x,y)=y^{7}\ast ^{y}x^{4}y^{5}=\frac{1}{10296}x^{4}y^{13}.
\label{eq:tr883}
\end{equation}

On using Eqs (\ref{eq:tr4.4}), (\ref{eq:tr881}), (\ref{eq:tr882}) and (\ref%
{eq:tr883}) we have
\begin{equation}
D_{2}(x,y)=-\frac{1}{101930400}x^{6}y^{24}  \label{eq:tr884}
\end{equation}%
We can easily see from Eq(\ref{eq:tr884}) that Eq (\ref{eq:tr880})
is an elliptic equation for all $(x_{,}y_{0}).$\newline Now, let
us extend the above results form single convolution to double
convolution as:

\noindent If we multiply the Eq (\ref{eq:tr4.2}) by polynomial
with double convolution
as $k(x,y)\ast \ast $ where $\displaystyle k(x,y)=\sum_{j=1}^{n}%
\sum_{i=1}^{m}x^{i}y^{j}$ then Eq (\ref{eq:tr4.2}) becomes
\begin{equation*}
k(x,y)\ast \ast \left[
a(x,y)u_{xx}+b(x,y)u_{xy}+c(x,y)u_{yy}+F(x,y,u,u_{x},u_{y})\right] =0
\end{equation*}%
\noindent where the symbol \ $\ast \ast $ \ indicates double convolution and
we shall classify more general form Eq (\ref{eq:tr770}) instead of Eq (\ref%
{eq:tr4.2}) by considering the function
\begin{equation}
D(x,y)=\left( k(x,y)\ast \ast b(x,y)\right) ^{2}-\left( k(x,y)\ast \ast
a(x,y)\right) \left( k(x,y)\ast \ast c(x,y)\right)  \label{eq:tr4.4}
\end{equation}%
\noindent First of all, we compute the coefficients of Eq
(\ref{eq:tr4.4}) as follow
\begin{equation*}
A(x,y)=k(x,y)\ast \ast a(x,y)=\sum_{j=1}^{n}\sum_{i=1}^{m}x^{i}y^{j}\ast
\ast \sum_{\beta =1}^{n}\sum_{\alpha =1}^{m}x^{\alpha }y^{\beta }
\end{equation*}%
\noindent by using double convolution definition and integral by part we
obtain%
\begin{equation}
A(x,y)=\int_{0}^{y}\sum_{j=1}^{n}\sum_{\beta =1}^{n}(y-\nu )^{j}\nu ^{\beta
}d\nu \times \int_{0}^{x}\sum_{i=1}^{m}\sum_{\alpha =1}^{m}(x-\mu )^{i}\mu
^{\alpha }d\mu  \label{eq:tr4.5}
\end{equation}%
\noindent the first integral in the right hand side of Eq (\ref{eq:tr4.5})
can be obtain as
\begin{equation}
\int_{0}^{y}\sum_{j=1}^{n}\sum_{\beta =1}^{n}(y-\nu )^{j}\nu ^{\beta }d\nu
=\sum_{j=1}^{n}\sum_{\beta =1}^{n}\frac{j!y^{\beta +j+1}}{(\beta +1)((\beta
+2)...(\beta +j+1)}  \label{eq:tr4.6}
\end{equation}%
\noindent in a similar way the second integral in the right hand side of Eq (%
\ref{eq:tr4.5}) is given by
\begin{equation}
\int_{0}^{x}\sum_{i=1}^{m}\sum_{\alpha =1}^{m}(x-\mu )^{i\text{ }}\mu
^{\alpha }d\mu =\sum_{i=1}^{m}\sum_{\alpha =1}^{m}\frac{i!x^{\alpha +i+1}}{%
(\alpha +1)((\alpha +2){...}(\alpha +i+1)}  \label{eq:tr4.7}
\end{equation}%
\noindent thus from Eqs (\ref{eq:tr4.6}) and (\ref{eq:tr4.7}) we
obtain the first coefficient of Eq (\ref{eq:tr4.4}) in the form of
\begin{equation*}
A(x,y)=\sum_{j=1}^{n}\sum_{\beta =1}^{n}\sum_{i=1}^{m}\sum_{\alpha =1}^{m}%
\frac{i!j!x^{\alpha +i+1}y^{\beta +j+1}}{\left( (\alpha +1)((\alpha +2){...}%
(\alpha +i+1)\right) \left( (\beta +1)((\beta +2){...}(\beta +j+1)\right) }
\end{equation*}%
\noindent\ similarly, for\ the coefficients of the second part in Eq (\ref%
{eq:tr4.4}) we have
\begin{equation*}
B(x,y)=\sum_{j=1}^{n}\sum_{l=1}^{n}\sum_{k=1}^{m}\sum_{i=1}^{m}\frac{%
i!j!x^{k+i+1}y^{l+j+1}}{\left( (k+1)((k+2){...}(k+i+1)\right) \left(
(l+1)((l+2){...}(l+j+1)\right) }
\end{equation*}%
\noindent and similarly,  the last coefficient of Eq
(\ref{eq:tr4.4}) given by
\begin{equation*}
C(x,y)=\sum_{j=1}^{n}\sum_{l=1}^{n}\sum_{k=1}^{m}\sum_{i=1}^{m}\frac{%
i!j!x^{\zeta +i+1}y^{\eta +j+1}}{\left( (k+1)((k+2){...}(k+i+1)\right)
\left( (l+1)((l+2){...}(l+j+1)\right) }
\end{equation*}%
then one can easily set up
\begin{equation}
D(x,y)=B^{2}-A(x,y)C(x,y).  \label{eq:tr4.8}
\end{equation}%
From Eq (\ref{eq:tr4.8}), one can see that if $D$ is positive then Eq (\ref%
{eq:tr4.8}) it called Hyperbolic, if $D$ is negative then Eq (\ref{eq:tr4.8}%
) it called Elliptic, otherwise parabolic. \newline
\newpage
\textbf{Generalized the Classification of Hyperbolic and Elliptic equations}%
\bigskip

\noindent In the this section we are going to generalized the classification
of hyperbaric and elliptic equations as follow. Firstly, we consider the one
dimensional wave equation with polynomial as follows
\begin{equation}
a(x,t)u_{tt}-c(x,t)u_{xx}=f_{1}(x,t)\ast \ast f_{2}(x,t)  \label{eq:tr4}
\end{equation}%
where $a,b$ are polynomials defined by $\displaystyle a(x,t)=\sum_{\beta
=1}^{n}\sum_{\alpha =1}^{m}x^{\alpha }t^{\beta }$ and $\displaystyle %
c(x,t)=\sum_{l=1}^{n}\sum_{k=1}^{m}x^{\mu }t^{\nu }$, now, we generalize Eq (%
\ref{eq:tr4}) by using finite double convolutional product as follows \
\begin{equation}
\dprod\limits_{r=1}^{l}p_{r}(x,t)\ast \ast
(a(x,t)u_{tt}-c(x,t)u_{xx})=f_{1}(x,t)\ast \ast f_{2}(x,t)  \label{eq:tr35}
\end{equation}%
where $r=1,2,{\ldots }l,$ and
\begin{equation*}
\displaystyle{\dprod\limits_{r=1}^{l}\left( p_{r}(x,t)\right)
=p_{1}(x,t)\ast \ast p_{2}(x,t)\ast \ast p_{3}(x,t)\ast \ast \ldots \ast
\ast p_{l}(x,t)}
\end{equation*}%
where $p_{1},p_{2},p_{3}...p_{l}$ polynomials defined by
\begin{equation}
p_{1}(x,t)=\sum_{j_{1}=1}^{n}\sum_{i_{1}=1}^{m}x^{i_{1}}t^{j_{1}},\text{ }%
p_{2}(x,t)=\sum_{j_{2}=1}^{n}\sum_{i_{2}=1}^{m}x^{i_{2}}t^{j_{2}},{\ldots },%
\text{ }p_{l}(x,t)=\sum_{j_{l}=1}^{n}\sum_{i_{l}=1}^{m}x^{i_{l}}t^{j_{l}}
\label{eq:tr36}
\end{equation}%
To compute the first coefficients of Eq (\ref{eq:tr35}), we start by $l=1,$
we obtain the coefficient in the form
\begin{equation*}
A_{1}(x,t)=p_{1}(x,t)\ast \ast
a(x,t)=\sum_{j_{1}=1}^{n}\sum_{i_{1}=1}^{m}x^{i_{1}}t^{j_{1}}\ast \ast
\sum_{\beta =1}^{n}\sum_{\alpha =1}^{m}x^{\alpha }t^{\beta }
\end{equation*}%
\noindent by using double convolution definition and integral by part we
obtain%
\begin{equation}
A_{1}(x,t)=\int_{0}^{t}\sum_{j_{1}=1}^{n}\sum_{\beta =1}^{n}(t-\zeta
)^{j_{1}}\zeta ^{\beta }d\zeta \times
\int_{0}^{x}\sum_{i_{1}=1}^{m}\sum_{\alpha =1}^{m}(x-\eta )^{i_{1}}\eta
^{\alpha }d\eta  \label{eq:tr4.55}
\end{equation}%
\noindent the first integral in the right hand side of Eq (\ref{eq:tr4.55})
can be obtain as
\begin{equation}
\int_{0}^{t}\sum_{j_{1}=1}^{n}\sum_{\beta =1}^{n}(t-\zeta )^{j_{1}}\zeta
^{\beta }d\zeta =\sum_{j_{1}=1}^{n}\sum_{\beta =1}^{n}\frac{j_{1}!t^{\beta
+j+1}}{(\beta +1)((\beta +2)...(\beta +j_{1}+1)}  \label{eq:tr4.66}
\end{equation}%
\noindent in a similar way the second integral in the right hand side of Eq (%
\ref{eq:tr4.55}) is given by
\begin{equation}
\int_{0}^{x}\sum_{i_{1}=1}^{m}\sum_{\alpha =1}^{m}(x-\eta )^{i_{1}}\eta
^{\alpha }d\eta =\sum_{i_{1}=1}^{m}\sum_{\alpha =1}^{m}\frac{i_{1}!x^{\alpha
+i+1}}{(\alpha +1)((\alpha +2){...}(\alpha +i_{1}+1)}  \label{eq:tr4.77}
\end{equation}%
\noindent thus from Eqs (\ref{eq:tr4.66}) and (\ref{eq:tr4.77}) we obtain
\begin{equation*}
A_{1}(x,t)=\sum_{j_{1}=1}^{n}\sum_{\beta
=1}^{n}\sum_{i_{1}=1}^{m}\sum_{\alpha =1}^{m}x^{\alpha +i_{1}+1}t^{\beta
+j_{1}+1},
\end{equation*}

where
\begin{equation*}
Q=\frac{i_{1}!j_{1}!}{\left( (\alpha +1)((\alpha +2){...}(\alpha
+i_{1}+1)\right) \left( (\beta +1)((\beta +2){...}(\beta +j_{1}+1)\right) }
\end{equation*}%
similarly, if we let $l=2,$ we have
\begin{equation*}
A_{2}(x,t)=A_{1}(x,t)\ast \ast p_{1}(x,t)=\sum_{j_{1}=1}^{n}\sum_{\beta
=1}^{n}\sum_{i_{1}=1}^{m}\sum_{\alpha =1}^{m}Qx^{\alpha +i_{1}+1}t^{\beta
+j_{1}+1}\ast \ast \sum_{j_{2}=1}^{n}\sum_{i_{2}=1}^{m}x^{i_{2}}t^{j_{2}}
\end{equation*}

similarly, as above we obtain the coefficient in case $l=2,$ as follows
\begin{equation*}
A_{2}(x,t)=\sum_{j_{1}=1}^{n}\sum_{j_{2}=1}^{n}\sum_{\beta
=1}^{n}\sum_{i_{1}=1}^{m}\sum_{i_{2}=1}^{m}\sum_{\alpha =1}^{m}\frac{K}{GL}%
x^{\alpha +i_{1}+i_{2}+1}t^{\beta +j_{1}+j_{2}+1}
\end{equation*}

where
\begin{eqnarray*}
K &=&i_{1}!i_{2}!j_{1}!j_{2}!\text{, }G=\left( (\alpha +1)(\alpha
+2)...(\alpha +i_{1}+1)\right) \left( \left( i_{1}+\alpha +2\right) \left(
i_{1}+\alpha +3\right) ...\left( i_{1}+\alpha +1+i_{2}\right) \right) \\
\text{and }L &=&\left( (\beta +1)(\beta +2)...(\beta +j_{1}+1)\right) \left(
((j_{1}+\beta +2))((j_{1}+\beta +3))...((j_{1}+\beta +1+j_{2}))\right) .
\end{eqnarray*}

By the same way we compute the coefficient in case $l=3$%
\begin{equation*}
A_{3}(x,t)=A_{2}(x,t)\ast \ast
p_{3}(x,t)=\sum_{j_{1}=1}^{n}\sum_{j_{2}=1}^{n}\sum_{\beta
=1}^{n}\sum_{i_{1}=1}^{m}\sum_{i_{2}=1}^{m}\sum_{\alpha =1}^{m}\frac{K}{GL}%
x^{\alpha +i_{1}+i_{2}+1}t^{\beta +j_{1}+j_{2}+1}\ast \ast
\sum_{j_{3}=1}^{n}\sum_{i_{3}=1}^{m}x^{i_{3}}t^{j_{3}}
\end{equation*}%
by using double convolution definition and integral by part we get%
\begin{equation*}
A_{3}(x,t)=\sum_{j_{1}=1}^{n}\sum_{j_{2}=1}^{n}\sum_{j_{3}=1}^{n}\sum_{\beta
=1}^{n}\sum_{i_{1}=1}^{m}\sum_{i_{2}=1}^{m}\sum_{i_{3}=1}^{m}\sum_{\alpha
=1}^{m}\frac{B}{MN}x^{\alpha +i_{1}+i_{2}+i_{3}+1}t^{\beta
+j_{1}+j_{2}+j_{3}+1},
\end{equation*}%
where
\begin{equation*}
B=i_{1}!i_{2}!i_{3}!j_{1}!j_{2}!j_{3}!,
\end{equation*}%
\begin{eqnarray*}
M &=&\left( (\alpha +1)(\alpha +2)...(\alpha +i_{1}+1)\right) \left( \left(
i_{1}+\alpha +2\right) \left( i_{1}+\alpha +3\right) ...\left( i_{1}+\alpha
+1+i_{2}\right) \right) \\
&&\times (\left( \left( i_{1}+i_{2}+\alpha +3\right) \left(
i_{1}+i_{2}+\alpha +4\right) ...\left( i_{1}+i_{2}+\alpha +1+i_{3}\right)
\right) )
\end{eqnarray*}

and
\begin{eqnarray*}
N &=&\left( (\beta +1)(\beta +2)...(\beta +j_{1}+1)\right) \left(
((j_{1}+\beta +2))((j_{1}+\beta +3))...((j_{1}+\beta +1+j_{2}))\right) \\
&&\times ((j_{1}+j_{2}+\beta +3))...((j_{1}+j_{2}+\beta +1+j_{3}))
\end{eqnarray*}

\noindent In general form the first coefficients of Eq (\ref{eq:tr35}) given
by%
\begin{eqnarray}
A(x,t) &=&\dprod\limits_{r=1}^{l}p_{r}(x,t)\ast \ast (a(x,t)  \notag \\
&=&\sum_{j_{1}=1}^{n}\sum_{i_{1}=1}^{m}{\ldots }\sum_{j_{1}=1}^{n}%
\sum_{i_{1}=1}^{m}\left( \sum_{\beta =1}^{n}\sum_{\alpha =1}^{m}\frac{\Phi
\Psi }{\Theta \Xi }x^{i_{1}+i_{2}+..+i_{l}+\alpha
+1}t^{j_{1}+j_{2}+...+j_{l}+\beta +1}\right) \ \   \label{eq:tr37}
\end{eqnarray}%
where
\begin{eqnarray*}
\Phi &=&i_{1}!i_{2}!...i_{l}!,\text{ \ }\Psi =j_{1}!j_{2}!...j_{l}!, \\
\Theta &=&\left( (\alpha +1)(\alpha +2)...(\alpha +i_{1}+1)\right) \left(
\left( i_{1}+\alpha +2\right) \left( i_{1}+\alpha +3\right) ...\left(
i_{1}+\alpha +1+i_{2}\right) \right) ... \\
&&...\left( (i_{1}+...+i_{l-1}+\alpha +1)...(i_{1}+...+i_{l}+\alpha
+1)\right)
\end{eqnarray*}%
and
\begin{eqnarray*}
\Xi &=&\left( (\beta +1)(\beta +2)...(\beta +j_{1}+1)\right) \left(
((j_{1}+\beta +2))((j_{1}+\beta +3))...((j_{1}+\beta +1+j_{2}))\right) ... \\
&&...\left( (j_{1}+j_{2}+...+j_{l-1}+\beta
+1)...(j_{1}+j_{2}+...+j_{l}+\beta +1)\right)
\end{eqnarray*}%
Similarly, the second coefficients of Eq (\ref{eq:tr35}) given by
\begin{eqnarray}
C(x,t) &=&\dprod\limits_{r=1}^{l}p_{r}(x,t)\ast \ast (c(x,t)  \notag \\
&=&\sum_{j_{1}=1}^{n}\sum_{i_{1}=1}^{m}{\ldots }\sum_{j_{1}=1}^{n}%
\sum_{i_{1}=1}^{m}\left( \sum_{\mu =1}^{n}\sum_{\nu =1}^{m}\frac{\digamma
\Upsilon }{\Delta \Lambda }x^{i_{1}+i_{2}+...+i_{l}+\mu
+1}t^{j_{1}+j_{2}+...+j_{l}+\nu +1}\right) \ \   \label{eq:tr38}
\end{eqnarray}%
where
\begin{eqnarray*}
\digamma &=&i_{1}!i_{2}!...i_{l}!,\text{ \ }\Upsilon =j_{1}!j_{2}!...j_{l}!,
\\
\Delta &=&\left( (\mu +1)(\mu +2)...(\mu +i_{1}+1)\right) \left( (\mu
+i_{1}+2)(\mu +i_{1}+3)...(\mu +i_{1}+i_{2}+1)\right) ... \\
&&...\left( (i_{1}+i_{2}+...+i_{l-1}+\mu +1)...(i_{1}+i_{2}+...+i_{l}+\mu
+1)\right)
\end{eqnarray*}%
and
\begin{eqnarray*}
\Lambda &=&\left( (\nu +1)(\nu +2)...(\nu +j_{1}+1)\right) \left( (\nu
+j_{1}+2)(\nu +j_{1}+3)...(\nu +j_{1}+j_{2}+1)\right) ... \\
&&...\left( (j_{1}+j_{2}+...+j_{l-1}+\nu +1)...(j_{1}+j_{2}+...+j_{l}+\nu
+1)\right)
\end{eqnarray*}%
Now, we going back to Eq(\ref{eq:tr4.8})
\begin{equation*}
D(x,t)=-A(x,t)C(x,t).
\end{equation*}%
Then we obtain $D(x,t)$ as follow
\begin{eqnarray}
D(x,t) &=&\sum_{j_{1}=1}^{n}\sum_{i_{1}=1}^{m}{\ldots }\sum_{j_{1}=1}^{n}%
\sum_{i_{1}=1}^{m}\left( \sum_{\beta =1}^{n}\sum_{\alpha =1}^{m}\frac{\Phi
\Psi }{\Theta \Xi }x^{i_{1}+i_{2}+..+i_{n}+\alpha
+1}t^{j_{1}+j_{2}+...+j_{m}+\beta +1}\right) \times  \notag \\
&&\sum_{j_{1}=1}^{n}\sum_{i_{1}=1}^{m}{\ldots }\sum_{j_{1}=1}^{n}%
\sum_{i_{1}=1}^{m}\left( \sum_{\mu =1}^{n}\sum_{\nu =1}^{m}\frac{\digamma
\Upsilon }{\Delta \Lambda }x^{i_{1}+i_{2}+..+i_{n}+\mu
+1}t^{j_{1}+j_{2}+...+j_{m}+\nu +1}\right)  \label{eq:tr40}
\end{eqnarray}%
We assume that all the coefficients $A(x,t),$ $C(x,t)$ are convergent. Now,
we can consider some particular cases: \newline

\noindent \textbf{Case(1):} If $i_{1}+i_{2}+...+i_{l}+j_{1}+j_{2}+...+j_{l},%
\alpha +\beta $ and \ $\mu +\nu $\ are odd, we classify Eq (\ref{eq:tr35})
by using Eq (\ref{eq:tr40}), \noindent since $\alpha +\beta $ and \ $\mu
+\nu $\ are odd, then the power of $x,t$ in Eq (\ref{eq:tr40}) is even, thus
we see that $D>0$ for all $(x_{0,}t_{0})\in
\mathbb{R}
^{2},$ \noindent Eq (\ref{eq:tr35}) Hyperbolic equation. Then for example if
we consider a particular non-constant coefficient wave equation in one
dimension in the form%
\begin{equation}
x^{5}t^{2}\ast \ast \left( x^{4}t^{3}\ast \ast x^{2}t^{7}\right)
u_{tt}-x^{5}t^{2}\ast \ast \left( x^{4}t^{3}\ast \ast x^{6}t^{5}\right)
u_{xx}=f(x,t)\ast \ast g(x,t)  \label{eq:tr40.}
\end{equation}

\noindent We follow the similar technique that was introduced by A. K\i l\i
\c{c}man and H. Eltayeb, see \cite{Adem}, then in particular we compute the
coefficients of Eq(\ref{eq:tr40.}) by using MAPLE 11 as
\begin{equation*}
A=x^{5}t^{2}\ast \ast x^{4}t^{3}\ast \ast x^{2}t^{7}=\frac{1}{1558311955200}%
x^{13}t^{14},
\end{equation*}%
\noindent where
\begin{eqnarray*}
\Phi  &=&i_{1}!i_{2}!=4!5! \\
\Psi  &=&j_{1}!j_{2}!=3!2! \\
\Theta  &=&\left( (\alpha +1)(\alpha +2)...(\alpha +i_{1}+1)\right) \left(
\left( i_{1}+\alpha +2\right) \left( i_{1}+\alpha +3\right) ...\left(
i_{1}+\alpha +2+i_{2}\right) \right)
\end{eqnarray*}%
and $\displaystyle\Xi =\left( (\beta +1)(\beta +2)...(\beta +j_{1}+1)\right)
\left( ((j_{1}+\beta +2))((j_{1}+\beta +3))...((j_{1}+\beta
+2+j_{2}))\right) ,$ where $i_{1}=4,$ $\alpha =2,$ $i_{2}=5,$ $j_{1}=3,$ $%
\beta =7$ and $j_{2}=2,$ thus%
\begin{equation*}
\frac{\Phi \Psi }{\Theta \Xi }=\frac{4!5!3!2!}{%
(3.4.5.6.7)(8.9.10.11.12.13)(8.9.10.11)(12.13.14)}=\frac{1}{1558311955200}
\end{equation*}%
\begin{equation*}
C=-x^{5}t^{2}\ast \ast x^{4}t^{3}\ast \ast x^{6}t^{5}=-\frac{1}{%
57058191590400}x^{17}t^{12},
\end{equation*}%
similarly, the factor of second coefficients is given by
\begin{equation*}
\frac{\digamma \Upsilon }{\Delta \Lambda }=\frac{1}{57058191590400},
\end{equation*}%
then%
\begin{equation}
D=-AC=\frac{1}{88914462097412421550080000}x^{30}t^{26}  \label{eq:tr40..}
\end{equation}%
We can easy see from Eq(\ref{eq:tr40..}) that Eq(\ref{eq:tr40.}) is
hyperbolic for all $(x_{0,}t_{0})\in
\mathbb{R}
^{2}.$\newline

\noindent \textbf{Case(2):} If $i_{1}+i_{2}+...+i_{l}+j_{1}+j_{2}+...+j_{l},%
\alpha +\beta $ and \ $\mu +\nu $ are even in Eq (\ref{eq:tr35}), and assume
$b=0$,then the power of $x,t$ in equation Eq (\ref{eq:tr40}) is even thus it
follows that for all points $(x_{0},t_{0})$ in the domain $%
\mathbb{R}
^{2}$ Eq (\ref{eq:tr35}) Hyperbolic.\newline
Now, if we consider another particular non-constant coefficient wave
equation in one dimension in the form of%
\begin{equation}
xt^{9}\ast \ast \left( x^{5}t^{3}\ast \ast x^{3}t^{7}\right)
u_{tt}-xt^{9}\ast \ast \left( x^{5}t^{3}\ast \ast x^{7}t^{5}\right)
u_{xx}=f(x,t)\ast \ast g(x,t)  \label{eq:tr41.}
\end{equation}
same as above example, we make calculations and obtain%
\begin{equation}
D(x,t)=-AC=\frac{1}{259841930424676205263257600000}x^{26}t^{40}
\label{eq:tr41..}
\end{equation}
From Eq(\ref{eq:tr41..}) we see that Eq(\ref{eq:tr41.}) is hyperbolic for
all point in $%
\mathbb{R}
^{2}.$\newline

\noindent \textbf{Case (3)}: If $i_{1}+i_{2}+...+i_{l}+j_{1}+j_{2}+...+j_{l}$
are odd$,\alpha +\beta $ and \ $\mu +\nu $\ are even in Eq (\ref{eq:tr35}),
then we are going to classify Eq (\ref{eq:tr35}) by using Eq (\ref{eq:tr40}%
), in this case $b=0$ therefore we see that the power of $x,t$ in Eq (\ref%
{eq:tr40}) are even then it follows that for all point $(x_{0},t_{0})$ in
the domain $%
\mathbb{R}
^{2}$ Eq (\ref{eq:tr35}) is Hyperbolic. \ As before, we consider simple
non-constant wave equation in one dimensional form:
\begin{equation}
x^{3}t^{4}\ast \ast \left( x^{6}t^{5}\ast \ast x^{3}t^{7}\right)
u_{tt}-x^{3}t^{4}\ast \ast \left( x^{6}t^{5}\ast \ast x^{7}t^{5}\right)
u_{xx}=f(x,t)\ast \ast g(x,t)  \label{eq:tr50}
\end{equation}%
and as above we see that%
\begin{equation}
D(x,t)=\frac{1}{26409729190021827098490961920000}x^{32}t^{34}
\label{eq:tr51}
\end{equation}%
Thus Eq(\ref{eq:tr51}) is positive for all points $(x_{0},t_{0})$ $\in
\mathbb{R}
^{2};$ thus Eq (\ref{eq:tr50}) is hyperbolic.\newline

\noindent \textbf{Case(4):} If $i_{1}+i_{2}+...+i_{l}+j_{1}+j_{2}+...+j_{l}$
are even$,\alpha +\beta $ and \ $\mu +\nu $\ are odd in Eq (\ref{eq:tr35}),
it is easy to see the power of $x,t$ in Eq (\ref{eq:tr40}) are even , then\
for all point $(x_{0},t_{0})$ in the domain $%
\mathbb{R}
^{2}$ Eq (\ref{eq:tr35}) is Hyperbolic.\newline

\noindent Secondly, we aim to examine the classification of generalized
elliptic equation; we consider the equations with non-constant coefficients
in the form of
\begin{equation}
a(x,y)u_{xx}+2b(x,y)u_{xy}+c(x,y)u_{yy}=f_{1}(x,y)\ast \ast f_{2}(x,y)
\label{eq:tr26}
\end{equation}%
where $a,b,c$ and $k,$ the polynomials and $\displaystyle{%
\dprod\limits_{r=1}^{z}\left( k_{r}(x,y) \right) ,}$ the finite double
convolution then Eq (\ref{eq:tr26}) becomes \
\begin{equation}
\dprod\limits_{r=1}^{z}k_{r}(x,y)\ast \ast \left(
a(x,y)u_{xx}+2b(x,y)u_{xy}+c(x,y)u_{yy}\right) =f_{1}(x,y)\ast \ast
f_{2}(x,y)  \label{eq:tr41}
\end{equation}%
where $r=1,2,{\ldots }z,$ and
\begin{equation*}
\dprod\limits_{r=1}^{z}\left( k_{r}(x,y)\right) =k_{1}(x,y)\ast \ast
k_{2}(x,y)\ast \ast k_{3}(x,y)\ast \ast \ldots \ast \ast k_{z}(x,y)
\end{equation*}%
where
\begin{equation*}
k_{1}(x,y) = \sum_{j_{1}=1}^{n}\sum_{i_{1}=1}^{m}x^{i_{1}}y^{j_{1}},\ \ \
k_{2}(x,y) = \sum_{j_{2}=1}^{n}\sum_{i_{2}=1}^{m}x^{i_{2}}y^{j_{2}},\ \ldots
, \ k_{z}(x,y) = \sum_{j_{l}=1}^{n}\sum_{i_{l}=1}^{m}x^{i_{z}}y^{j_{z}}
\end{equation*}
Let us compute the coefficients of Eq (\ref{eq:tr41}), by similar way that
used to compute the coefficients of \ equation (\ref{eq:tr35}), the first
coefficient given by%
\begin{eqnarray*}
A(x,y) &=&\dprod\limits_{r=1}^{z}k_{r}(x,y)\ast \ast a(x,y) \\
&=&\sum_{j_{1}=1}^{n}\sum_{i_{1}=1}^{m}{\ldots }\sum_{j_{1}=1}^{n}%
\sum_{i_{1}=1}^{m}\left( \sum_{\beta =1}^{n}\sum_{\alpha =1}^{m}P_{\alpha
\beta }\frac{\Gamma \Omega }{MN}x^{i_{1}+i_{2}+..+i_{z}+\alpha
+1}y^{j_{1}+j_{2}+...+j_{z}+\beta +1}\right) \ \
\end{eqnarray*}%
where
\begin{eqnarray*}
\Gamma &=&i_{1}!i_{3}!...i_{2z-1}!,\text{ \ }\Omega
=j_{1}!j_{2}!...j_{2z-1}!, \\
M &=&\left( (\alpha +1)(\alpha +2)...(\alpha +i_{1}+1)\right) \left( (\alpha
+i_{1}+2)(\alpha +i_{1}+3)...(\alpha +i_{1}+i_{2}+1)\right) ... \\
&&...\left( (i_{1}+i_{2}+...+i_{z-1}+\alpha
+2)...(i_{1}+i_{2}+...+i_{z}+\alpha +1)\right)
\end{eqnarray*}%
and
\begin{eqnarray*}
N &=&\left( (\beta +1)(\beta +2)...(\beta +j_{1}+1)\right) \left( (\beta
+j_{1}+2)(\beta +j_{1}+3)...(\beta +j_{1}+j_{2}+1)\right) ... \\
&&...\left( (j_{1}+j_{2}+...+j_{z-1}+\beta
+2)...(j_{1}+j_{2}+...+j_{z}+\beta +1)\right)
\end{eqnarray*}%
Similarly, the second coefficient given by
\begin{eqnarray}
B(x,t) &=&\dprod\limits_{r=1}^{z}k_{r}(x,y)\ast \ast b(x,y)  \notag \\
&=&\sum_{j_{1}=1}^{n}\sum_{i_{1}=1}^{m}{\ldots }\sum_{j_{1}=1}^{n}%
\sum_{i_{1}=1}^{m}\left( \sum_{t=1}^{n}\sum_{s=1}^{m}\frac{FG}{QE}%
x^{i_{1}+i_{2}+..+i_{z}+s+1}y^{j_{1}+j_{2}+...+j_{z}+t+1}\right) \ \
\label{eq:tr44}
\end{eqnarray}%
where
\begin{eqnarray*}
F &=&i_{1}!i_{3}!...i_{2z-1}!,\text{ \ }G=j_{1}!j_{2}!...j_{2z-1}!t!, \\
Q &=&\left( (s+1)(s+2)...(s+i_{1}+1)\right) \left(
(s+i_{1}+2)(s+i_{1}+3)...(s+i_{1}+i_{2}+1)\right) ... \\
&&...\left(
(i_{1}+i_{2}+...+i_{z-1}+s+2)...(i_{1}+i_{2}+...+i_{z}+s+1)\right)
\end{eqnarray*}%
\begin{eqnarray*}
E &=&\left( (t+1)(t+2)...(t+j_{1}+1)\right) \left(
(t+j_{1}+2)(t+j_{1}+3)...(t+j_{1}+j_{2}+1)\right) ... \\
&&...\left(
(j_{1}+j_{2}+...+j_{z-1}+t+2)...(j_{1}+j_{2}+...+j_{z}+t+1)\right)
\end{eqnarray*}%
and $\displaystyle{s=\frac{\alpha +k}{2}}$, $\displaystyle{t=\frac{\beta +l}{%
2}.}$ Also the last coefficient given by%
\begin{eqnarray}
C(x,t) &=&\dprod\limits_{r=1}^{z}k_{r}(x,y)\ast \ast c(x,y)  \notag \\
&=&\sum_{j_{1}=1}^{n}\sum_{i_{1}=1}^{m}{\ldots }\sum_{j_{1}=1}^{n}%
\sum_{i_{1}=1}^{m}\left( \sum_{t=1}^{n}\sum_{s=1}^{m}q_{kl}\frac{FG}{QE}%
x^{i_{1}+i_{2}+..+i_{z}+k+1}y^{j_{1}+j_{2}+...+j_{z}+l+1}\right) \ \ \ \
\label{eq:tr45}
\end{eqnarray}%
where
\begin{eqnarray*}
F &=&i_{1}!i_{3}!...i_{2z-1}!,\text{ \ }G=j_{1}!j_{2}!...j_{2z-1}!l!, \\
Q &=&\left( (k+1)(k+2)...(k+i_{1}+1)\right) \left(
(k+i_{1}+2)(k+i_{1}+3)...(k+i_{1}+i_{2}+1)\right) ... \\
&&...\left(
(i_{1}+i_{2}+...+i_{z-1}+k+2)...(i_{1}+i_{2}+...+i_{z}+k+1)\right)
\end{eqnarray*}%
\bigskip and
\begin{eqnarray*}
E &=&\left( (l+1)(l+2)...(l+j_{1}+1)\right) \left(
(l+j_{1}+2)(l+j_{1}+3)...(l+j_{1}+j_{2}+1)\right) ... \\
&&...\left(
(j_{1}+j_{2}+...+j_{z-1}+l+2)...(j_{1}+j_{2}+...+j_{z}+l+1)\right)
\end{eqnarray*}%
We assume that all the coefficients $A(x,y),B(x,y)$ $C(x,y)$ are convergent.
Now we can easily check the Eq (\ref{eq:tr41}) whether it is elliptic or
not. Now we have the following two cases.\newline

\noindent \textbf{(i):} If $i_{1}+i_{2}+...+i_{z}+j_{1}+j_{2}+...+j_{z},s+t,%
\alpha +\beta $ and $k+l$ are odd in the Eq (\ref{eq:tr41}), the power $s=%
\frac{\alpha +k}{2}$ and $t=\frac{\beta +l}{2}$ and the power of $x$ and $y$
in polynomials $a(x,y),$ $c(x,y)$ either even or odd, also the coefficient
of two polynomials $a(x,y),$ $c(x,y)$ have the same sign, now we are going
to study the classification of Eq (\ref{eq:tr41}), now the power of $x,y$ in
$B(x,y)^{2}=$ the power of $A(x,y)C(x,y)$ and the coefficient of $%
A(x,y)C(x,y)>1$ then the power is even, then\ for all point $(x_{0},y_{0})$
in the domain $%
\mathbb{R}
^{2}$ \noindent the Eq (\ref{eq:tr41}) is elliptic equation. For example in
particular we can have
\begin{equation}
x^{7}y^{4}\ast \ast \left[ \left( x^{3}y^{2}\ast \ast x^{2}y^{3}\right)
u_{xx}+ \left( x^{3}y^{2}\ast \ast x^{3}y^{4}\right) u_{xy}+ \left(
x^{3}y^{2}\ast \ast x^{4}y^{5}\right) u_{yy}\right]=f(x,y)\ast \ast g(x,y)
\label{eq:tr31}
\end{equation}%
and compute the coefficient of Eq (\ref{eq:tr31}) by using Eq (\ref{eq:tr4.8}%
), Eq (\ref{eq:tr44}) and Eq (\ref{eq:tr45}), we obtain
\begin{equation}
D(x,y)=-\frac{23}{466987721099855155200000000}x^{30}t^{24}  \label{eq:tr32}
\end{equation}

\noindent Then it is easy to see that Eq(\ref{eq:tr32}) is always negative
for all $(x_{0},y_{0})$ in the domain $%
\mathbb{R}
^{2}$, and thus Eq (\ref{eq:tr31}) is an elliptic equation.\newline

\noindent \textbf{(ii):} If $i_{1}+i_{2}+...+i_{z}+j_{1}+j_{2}+...+j_{z},$ $%
s+t,\alpha +\beta $ and $k+l$ are even and Eq (\ref{eq:tr41}) to be elliptic
equation at the condition since the power $s=\frac{\alpha +k}{2}$ and $t=%
\frac{\beta +l}{2}$, also the coefficients of two polynomials $a(x,y),$ $%
c(x,y)$ have the same sign, thus similar to the previous section the
classification of equation the Eq\ (\ref{eq:tr41}) on using Eq(\ref{eq:tr4.8}%
) and since the power of $x$ and $y$ in Eq (\ref{eq:tr41}) are even, then
for all point $(x_{0},y_{0})$ in the domain $%
\mathbb{R}
^{2}$ \noindent Eq (\ref{eq:tr41}) is elliptic equation.\newline

\noindent \textbf{Conclusion}

\noindent We note that the classification of generalized hyperbolic and
elliptic equations with non-constant coefficients are still similar with
original equations, that is, they are invariant after convolutions product.

\end{document}